\documentclass[a4paper,11pt]{article} %

\usepackage[utf8]{inputenc}
\usepackage{amsfonts,amsmath}%
\usepackage{graphicx}
\usepackage{algorithm}

\newcommand{\naturaln}{\ensuremath{{\mathbb{N}}}}
\newcommand{\realr}{\ensuremath{{\mathbb{R}}}}
\DeclareMathOperator{\prox}{prox}

\DeclareMathOperator{\grad}{\nabla}

\newtheorem{predefinition}{Definition}
  \newenvironment{definition}%
   {\begin{predefinition}\upshape}{\end{predefinition}}

\newtheorem{pretheorem}{Theorem}
  \newenvironment{theorem}%
   {\begin{pretheorem}\upshape}{\end{pretheorem}}

\newtheorem{prelemma}{Lemma}
  \newenvironment{lemma}%
   {\begin{prelemma}\upshape}{\end{prelemma}}
   
\newcommand{\proof}{Proof:}
\newcommand{\qed}{\hfill $\Box$}

\begin{document}
\title{On starting and stopping criteria for nested primal-dual iterations}
\author{Jixin Chen \and Ignace Loris}
\date{}


\maketitle


\begin{abstract}
The importance of an adequate inner loop starting point (as opposed to a sufficient inner loop stopping rule) is discussed in the context of an numerical optimization algorithm consisting of nested primal-dual proximal-gradient iterations. 
While the number of inner iterations is fixed in advance, convergence of the whole algorithm is still guaranteed by virtue of a warm-start strategy for the inner loop, showing that inner loop ``starting rules" can be just as effective as ``stopping rules'' for guaranteeing convergence.
The algorithm itself is applicable to the numerical solution of convex optimization problems defined by the sum of a differentiable term and two possibly non-differentiable terms. One of the latter terms should take the form of the composition of a linear map and a proximable function, while the differentiable term needs an accessible gradient. The algorithm reduces to the classical proximal gradient algorithm in certain special cases and it also generalizes other existing algorithms. In addition, under some conditions of strong convexity, we show a linear rate of convergence.
\end{abstract}


\section{Introduction}

Iterative optimization algorithms are based on the availability of simple building blocks related to the cost function that needs minimization. These building blocks, such as e.g. gradients and hessians in case of smooth optimization, should be easy to compute as they are evaluated in every step of the iteration process. 

In constrained convex optimization projections onto convex sets play a crucial role, such as in the projected gradient algorithm \cite{Goldstein1964}:
\begin{equation}
u_{n+1}=P_C(u_n-\alpha \grad f(u_n))\qquad\text{for}\qquad \min_{u\in C}f(u).
\end{equation}
Here $f$ is a real-valued, convex, differentiable function and $C$ a non-empty closed convex set in $\realr^d$.
More generally, the proximal-gradient algorithm (see e.g. \cite{ComWa2005})
\begin{equation}\label{eq:proxgrad}
u_{n+1}=\prox_{\alpha h}(u_n-\alpha \grad f(u_n)),
\end{equation}
(with $u_0$ arbitrary and $\alpha>0$ a step length parameter) can be applied to the numerical solution of the optimization problem 
\begin{equation}\label{eq:problem1}
\min_{u\in \realr^d}f(u)+h(u)
\end{equation}
where $h$ is a convex proper lower semi-continuous function.
Apart from the gradient of the differentiable part, here one also needs the \emph{proximal operator} of the non-differentiable function $h$, which was introduced in \cite{Moreau1965} (see also Definition~\ref{def:prox}) and for which explicit (and easy to evaluate) expressions exist for several useful cases \cite{ComWa2005}. 
One convergence result among many states that the algorithm (\ref{eq:proxgrad}) will converge to a minimizer (if one exists) when $0<\alpha<2/L$ where $L$ is the Lipschitz constant of the gradient of $f$ \cite{ComWa2005}.

A large number of generalizations of the proximal gradient algorithm exist, such as e.g. versions with variable metrics \cite{Combettes2014a,Chouzenoux-etal-2014} that exchange the Euclidean distance in the definition of the proximal operator for other ones, depending on the iteration step $n$ (see also \cite{Condat2013,Combettes2014b}).

In this note we are interested in an optimization problem defined by a cost function which consists of three parts instead of two:
\begin{equation}\label{eq:problem2}
\min_{u\in\realr^d} f(u)+g(Au)+h(u),
\end{equation}
where $f$ and $h$ are as before and where $A$ is a linear map and $g$ is a convex proper lower semi-continuous function.
If the proximal operator of $\alpha h+\alpha g\circ A$ were available, then algorithm (\ref{eq:proxgrad}) could be used for its solution (replacing $\prox_{\alpha h}$ in (\ref{eq:proxgrad}) by $\prox_{\alpha h+\alpha g\circ A}$) as follows:
\begin{equation}\label{eq:proxgrad2}
u_{n+1}=\prox_{\alpha h+\alpha g\circ A}(u_n-\alpha \grad f(u_n)),
\end{equation}
(with $u_0$ arbitrary).
We will however assume that the proximal operator of $\alpha h+\alpha g\circ A$ is not explicitly available, rendering algorithm (\ref{eq:proxgrad2}) ineffective. Still, we will suppose that the proximal operators $\prox_{\alpha h}$, $\prox_{\alpha g}$ and the linear operator $A$ separately are at our disposal (in many cases of practical interest the latter proximal operators are easier to compute in closed form than the former one). Such a problem has been studied in \cite{Zhang.Burger.ea2011,Combettes2012,Condat2013,Chen2016}.

It is well-known that the proximal operator appearing in (\ref{eq:proxgrad2}) can itself be found using an iterative algorithm based on dual variables (see for e.g. \cite{Chambolle2004,Chambolle2005,Aujol2009} for a special case in the area of mathematical imaging). Hence, a nested algorithm (i.e. combining an inner loop with an outer loop) can be a straight-forward way of tackling the described problem. Such nested primal-dual algorithms have already been used in practice for solving large scale optimization problem in mathematical imaging and signal processing \cite{Beck.Teboulle2009}. 

Using an inner loop for the calculation of the proximal operator invariably introduces numerical error in the outer loop. In general, convergence of the proximal gradient algorithm (\ref{eq:proxgrad2}) is robust with respect to errors of the proximal operator, in as much as the sum of all errors is finite \cite{ComWa2005}. However, such a condition is hard to verify in practice. Other (verifiable) conditions have also been proposed \cite{Bonettini2015,Bonettini2016}. The effects of inexact computation on accelerated proximal algorithms have been studied in \cite{Salzo2012,Schmidt2011}. 

In this paper we fix the number of inner iterations in advance (thereby completely avoiding the need to check a sufficient inner loop termination condition, while potentially losing control on the accuracy of the approximation of the proximal operator), but use a feedback procedure to guarantee the overall convergence.
The main goal of this paper therefore is to provide a rigorous convergence analysis of a nested iterative algorithm under an a priori finite termination condition of the inner loop (i.e. number of inner iterations is fixed in advance) with inner loop starting point feedback.

A ``warm start'' strategy is any method which uses the numerical result of one optimization problem as a starting point for a different, but closely related or perturbed, one. Such strategies are often used (e.g. for computing the solutions of a whole parameter family of optimization problems such as in \cite{Hu2016}), but theoretical guarantees or results are lacking. E.g. a warm start strategy is used in an inner loop in the proximal gradient ordered subsets framework applied to computer tomography in \cite{Rose2015}, but this is only briefly mentioned in an accompanying technical paper \cite[Algorithm 4 and below]{Rose2016}. Thus nested algorithms have already been proposed in the context of proximal algorithms; however, the convergence analysis presented here is novel, and puts the use of such algorithms on a firmer footing, bypassing a need to rely e.g. on the summability of errors in intermediate computations.

In the following section we will write a specific nested primal-dual algorithm applicable to the described problem. It uses only gradients of the differentiable term, the linear map $A$ (and its transpose) and the proximal operators of $g$ and $h$. In addition to its convergence we also prove a geometric convergence rate (under the additional assumption of strong convexity).
The nested primal-dual algorithm could be interpreted as a generalization of the algorithm proposed in \cite{Loris.Verhoeven2011} and further developed in \cite{Chen2013,Chen2016} in the sense that it could be identified as corresponding to just a single inner iteration in the algorithm discussed below. In that case, there is no real inner ``loop'' and the issue of its starting and stopping rule is absent. In this weak sense, the present proof of convergence generalizes the ones found in \cite{Loris.Verhoeven2011,Chen2016}.

%
%
%
%

\section{Nested Primal-Dual Proximal Gradient Algorithm}\label{sec:algorithm}

In the remainder of the paper we assume that $f:\realr^{d}\rightarrow \realr$ is convex, differentiable and that the gradient of $f$ is Lipschitz continuous (constant $L$). We also assume that $h:\realr^{d}\rightarrow \realr\cup \{+\infty\}$ and $g:\realr^{d'}\rightarrow \realr\cup \{+\infty\}$ are proper, convex, lower semi-continuous functions. Finally $A:\realr^{d}\rightarrow \realr^{d'}$ is a linear map, and $\|A\|$ signifies the largest singular value of the matrix $A$.

We start by recalling a number of well-known definitions and properties which are necessary for the derivation of the algorithms and for proving their convergence.
\begin{definition} Let $h:\realr^d\rightarrow\realr\cup\{+\infty\}$ be a convex, proper, lower semi-continuous function. The subdifferential of $h$ at the point $u$ is defined as the set
\begin{equation}
\partial h(u)=\{w\in\realr^{d}\quad |\quad h(v)\geq h(u)+\langle w,v-u\rangle\quad \forall v\in\realr^{d}\}.
\end{equation}
\end{definition}
It is easy to see that $\hat u$ is a minimizer of $h$ if and only if $0\in\partial h(\hat u)$. Also, under mild conditions \cite{Bauschke2011}, one can show that $\partial (h_1+h_2)(u)=\partial h_1(u)+\partial h_2(u)$ and $\partial (g\circ A)(u)=A^T\partial g(Au)$ where $A$ is a linear map.

\begin{definition}\label{def:prox} Let $h:\realr^d\rightarrow\realr\cup\{+\infty\}$ be a convex, proper, lower semi-continuous function. The proximal operator of $h$  is defined as:
\begin{equation} \label{eq:defprox}
\prox_{h}(a)=\arg\min_{u\in\realr^d} \frac{1}{2}\|u-a\|_2^2+ h(u).
\end{equation}
\end{definition}
The proximal operator is a nonexpansive map (Lipschitz continuous with constant $1$) defined on all of $\realr^d$.

\begin{definition}\label{def:fenchel} Let $h:\realr^d\rightarrow\realr\cup\{+\infty\}$ be a convex, proper, lower semi-continuous function. The Fenchel dual of $h$ is is defined as:
\begin{equation}
h^\ast:\realr^d\rightarrow\realr\cup\{+\infty\}: h^\ast(w)=\sup_u\langle w,u\rangle-h(u).
\end{equation}
\end{definition}
It is again a convex proper lower-semicontinuous function. In fact, on this class, the Fenchel transform is its own inverse: $(h^\ast)^\ast=h$.

Next we present some of the classical results of \cite{Moreau1965} under the form of the following lemmas.
\begin{lemma}\label{theoremMoreau}\label{subdiffproxprop2}
Let $h:\realr^d\rightarrow\realr\cup\{+\infty\}$ be a convex, proper, lower semi-continuous function. The following are equivalent:
\begin{enumerate}
\item $u=\prox_{\alpha h}(u+\alpha w)$ for any $\alpha>0$
\item $w\in\partial h(u)$
\item $h(u)+h^\ast(w)=\langle w,u\rangle$ \label{eq:Moreau}
\item $u\in\partial h^\ast(w)$
\item $w=\prox_{\beta h^\ast}(\beta u+w)$ for any $\beta>0$
\end{enumerate}
Furthermore, proximal operators of primal and dual functions $h$ and $h^\ast$ are related by Moreau's decomposition:
\begin{displaymath}
\prox_{\alpha h} (u)+\alpha \prox_{\alpha^{-1}h^\ast}(\alpha^{-1}u)=u\qquad \forall u\in\realr^d
\end{displaymath}
and any $\alpha>0$. It therefore suffices to know $\prox_h(a)$ in order to compute $\prox_{h^\ast}(a)$ and vice-versa.
\end{lemma}
\proof See \cite{Moreau1965}.\qed

Finally, we will need some further results on the Moreau envelope of a function.
\begin{definition}\label{def:envelope}
Let $\alpha>0$ and $h:\realr^d\rightarrow\realr\cup\{+\infty\}$ be convex proper and lower semi-continuous. The Moreau envelope of $h$ (of index $\alpha$) is defined as :
\begin{displaymath}
\begin{array}{lcl}
\hat h_\alpha:\realr^d\rightarrow\realr:\hat h_\alpha(u)&=& \min_{v\in\realr^d}\frac{1}{2\alpha}\|v-u\|_2^2+h(v) \\[3mm]
&=&\frac{1}{2\alpha}\|\prox_{\alpha h}(u)-u\|_2^2+h(\prox_{\alpha h}(u)).
\end{array}
\end{displaymath}
\end{definition}

\begin{lemma}\label{lemma:envelope}
The Moreau envelope admits the following properties:
\begin{enumerate}
\item $\hat h_\alpha$ is convex, proper and lower semi-continuous;
\item $\hat h_\alpha$ is differentiable and $\grad \hat h_\alpha(u)=\alpha^{-1}(u-\prox_{\alpha h}(u))$ (Lipschitz, constant $\alpha^{-1}$);\label{point2}
\item $\hat h_\alpha(u)+\widehat{(h^\ast)}_{1/\alpha}(u/\alpha)=\frac{1}{2\alpha}\|u\|_2^2$\quad (hence $\frac{1}{2\alpha}\|u\|_2^2-\hat h_\alpha(u)$ is a convex function).\label{point1}
\end{enumerate}
\end{lemma}
\proof See \cite{Moreau1965}.\qed

Our goal for the remainder of this section is to write an approximate version of the proximal gradient algorithm (\ref{eq:proxgrad2}), for problem (\ref{eq:problem2}). Informally, it takes the form:
\begin{equation}\label{eq:approxalgorithm}
\begin{array}{lcl}
u_{n+1}\approx \prox_{\alpha h+\alpha g\circ A} (u_n-\alpha \grad f(u_n)).
\end{array}
\end{equation}
In particular, we aim to approximate the proximal operator $\prox_{\alpha h+\alpha g\circ A} $ using an iterative calculation using dual variables. This calculation will involve the proximal operators of $\alpha h$ and of $\beta\alpha^{-1}g^\ast$, and the linear map $A$.

\begin{lemma}\label{lemma:dual}
The proximal operator of $\alpha h+ \alpha g\circ A$ evaluated at some point $a$, defined as:
\begin{displaymath} 
\hat a=\prox_{\alpha h+\alpha g\circ A} (a)=\arg\min_{u\in\realr^{d}} \frac{1}{2}\|u-a\|_2^2+\alpha h(u)+ \alpha g(Au),
\end{displaymath}
can be computed as $\hat a=\prox_{\alpha h}(a- \alpha A^T \hat v)$ where $\hat v$ is the limit of the sequence $(v^k)_{k\in\naturaln}$ defined by the iteration
\begin{equation}\label{eq:dualalg}
v^{k+1}=\prox_{\beta\alpha^{-1} g^\ast}(v^k+\beta\alpha^{-1} A\prox_{\alpha h}(a- \alpha A^T v^k)),
\end{equation}
for step size $0<\beta<2/\|A\|^2$ and arbitrary $v^0$.
\end{lemma}
\proof
Writing out the variational equations that determine the minimizer $\hat a$, one finds:
\begin{displaymath}
\hat a-a+ \alpha \hat w+ \alpha A^T \hat v=0\qquad\text{with}\quad \hat w\in\partial h(\hat a)\quad\text{and}\quad \hat v\in\partial g(A\hat a).
\end{displaymath}
Using Lemma~\ref{theoremMoreau} the latter inclusions can equivalently be written as a (non-linear) equations:
\begin{displaymath}
\hat a=\prox_{\alpha h}(\hat a+\alpha\hat w)\qquad\hat v=\prox_{\beta\alpha^{-1} g^\ast}(\hat v+\beta\alpha^{-1} A\hat a)
\end{displaymath}
where $\beta>0$ is an arbitrary parameter and $g^\ast$ is the Fenchel dual of $g$.
In other words $\hat a$ and $\hat v$ are determined by the equations:
\begin{displaymath}
\hat a=\prox_{\alpha h}(a- \alpha A^T \hat v)\qquad\text{and}\qquad \hat v=\prox_{\beta\alpha^{-1} g^\ast}(\hat v+\beta\alpha^{-1} A\hat a)
\end{displaymath}
for some $\beta>0$. Now the variable $\hat a$ can be eliminated from the second equation, yielding:
\begin{displaymath}
\hat v=\prox_{\beta\alpha^{-1} g^\ast}(\hat v+\beta\alpha^{-1} A\prox_{\alpha h}(a- \alpha A^T \hat v)).
\end{displaymath}

This equation for $\hat v$ can be solved using the fixed point iteration (\ref{eq:dualalg}). Indeed, if we set 
\begin{displaymath}
\varphi(u)=\frac{1}{2\alpha}\|u\|_2^2-\hat h_\alpha (u)
\end{displaymath}
(a convex differentiable function according to Lemma~\ref{lemma:envelope}, point~\ref{point1}) and 
\begin{displaymath}
\psi(v)=\alpha^{-1}\varphi(a-\alpha A^T v)
\end{displaymath}
(also a convex differentiable function), we see that the gradient of $\psi$ is (Lemma~\ref{lemma:envelope}, point~\ref{point2}):
\begin{displaymath}
\grad \psi(v)=\alpha^{-1} (-\alpha A)\grad \varphi(a-\alpha A^T v)=-\alpha^{-1} A \prox_{\alpha h}(a-\alpha A^T v).
\end{displaymath}
A Lipschitz constant of the gradient of $\psi$ is $\|A\|^2$.
Hence iteration (\ref{eq:dualalg}) is just the proximal gradient algorithm (\ref{eq:proxgrad}) applied to the ``dual problem'':
\begin{equation}\label{eq:dualproblem}
\min_v \psi(v)+\alpha^{-1}g^\ast(v)
\end{equation}
and therefore converges for $0<\beta<2/\|A\|^2$.
\qed

Introducing further auxiliary variables $u^k$ it is also possible to write iteration (\ref{eq:dualalg}) as:
\begin{equation}\label{dualalg1}
\begin{array}{l}
\text{for\ }k:0,1\ldots\\
\left\{
\begin{array}{l}
u^k=\prox_{\alpha h}(a- \alpha A^T v^k)\\
v^{k+1}=\prox_{\beta\alpha^{-1} g^\ast}(v^k+\beta\alpha^{-1} Au^k)\qquad 
\end{array}
\right.
\end{array}
\end{equation}
for step size $0<\beta<2/\|A\|^2$ and arbitrary $v^0$. As the sequence $(u^k)_{k\in\naturaln}$ in (\ref{dualalg1}) converges to $\prox_{\alpha h+\alpha g\circ A}(a)$, so does the sequence of averages. One can therefore write the following algorithm for approximating $\prox_{\alpha h+\alpha g\circ A}(a)$:
\begin{equation}\label{eq:dualalg2}
\begin{array}{l}\text{for\ }k:0\ldots k_{\max}-1\\
\left\{
\begin{array}{l}
u^k=\prox_{\alpha h}(a- \alpha A^T v^k)\\
v^{k+1}=\prox_{\beta\alpha^{-1} g^\ast}(v^k+\beta\alpha^{-1} Au^k)\qquad 
\end{array}
\right.\\[5mm]
u^{k_{\max}}=\prox_{\alpha h}(a- \alpha A^T v^{k_{\max}})\\[2mm]
\prox_{\alpha h+\alpha g\circ A}(a)\approx \sum_{k=1}^{k_{\max}}u^k/k_{\max}
\end{array}
\end{equation}
for some choice of $v^0$ and $k_{\max}\in\naturaln$.

Instead of imposing an implicit stopping rule on the iteration (\ref{eq:dualalg2}), such as e.g. requiring that $\|v^{k+1}-v^k\|_2<\epsilon$, we opt to fix the number of iterations $k_{\max}$ in advance. In general, this means that there is no guarantee as to the quality of the approximation (\ref{eq:dualalg2}). Indeed, the starting point could be chosen unfavorably. 

If $A, A^T$, $\prox_{\alpha h}$ and $\prox_{\beta\alpha^{-1} g^\ast}$ are available, algorithm (\ref{eq:dualalg2}) can be used to compute (an approximation of) the proximal operator present in algorithm (\ref{eq:approxalgorithm}). By replacing $a$ in (\ref{eq:dualalg2}) by $u_n-\alpha \grad f(u_n)$ we arrive at 
Algorithm~\ref{alg:npdaaverage}. We will systematically use subscripted $n$ as outer iteration index, and superscripted $k$ as inner iteration index.

\begin{algorithm} 
\caption{Nested primal dual algorithm}\label{alg:npdaaverage}
Choose $u_0, v_0^0$, $0<\alpha<2/L$, $0<\beta< 1/\|A\|^2$, $k_{\max}\in\naturaln_0$.
\begin{equation}\label{eq:alg4}
\begin{array}{l}
\text{for\ }n: 0,1,\ldots:\\[3mm]
\left\{
\begin{array}{l}
\text{for\ }k:0\ldots k_{\max}-1:\\[2mm]
\left\{
\begin{array}{l}
v_{n}^0=v_{n-1}^{k_{\max}}\qquad\qquad\text{for}\ n>0\\[3mm]
u_{n}^k=\prox_{\alpha h}(u_n-\alpha \grad f(u_n)-\alpha A^T v_{n}^k)\\[3mm]
v_{n}^{k+1}=\prox_{\beta\alpha^{-1} g^\ast}\left( v_{n}^{k}+\beta\alpha^{-1}Au_{n}^k\right) 
\qquad 
\end{array}
\right.\\[13mm]
u_{n}^{k_{\max}}=\prox_{\alpha h}(u_n-\alpha \grad f(u_n)-\alpha A^T v_{n}^{k_{\max}})\\[3mm]
u_{n+1}=\sum_{k=1}^{k_{\max}} u_{n}^{k}/k_{\max}
\end{array}
\right.
\end{array}
\end{equation}
\end{algorithm} 

It is important to note that, in the proposed nested algorithm, the inner loop starts with the outcome of the previous inner loop:  $v_{n}^0=v_{n-1}^{k_{\max}}$, and that the number of inner iterations $k_{\max}$ is fixed in advance. It is the former choice, rather than ``sufficient" inner iterations, that will allow use to prove convergence of this nested algorithm.


We remark that all iterates $u_n$ ($n\geq 1$) in Algorithm~\ref{alg:npdaaverage} are in the domain of $h$, but not necessarily in the domain of $g\circ A$. In the special case $h=0$ and $k_{\max}=1$ (just one inner iteration) Algorithm~\ref{alg:npdaaverage} reduces to
\begin{equation}\label{eq:alg3}
\left\{
\begin{array}{l}
z_n=u_n-\alpha \grad f(u_n)\\[3mm]
v_{n+1}=\prox_{\beta\alpha^{-1} g^\ast}\left( v_n+\beta\alpha^{-1}A(z_n-\alpha A^Tv_n)\right) \\[3mm]
u_{n+1}=z_n-\alpha A^Tv_{n+1},
\end{array}
\right.
\end{equation}
which was proposed in \cite{Loris.Verhoeven2011} and further studied in \cite{Chen2013,Chen2016}. It was also interpreted in \cite{Combettes2014} (see also \cite{Condat2013,Condat2014}) as a special case of a novel scheme extending several classical ones, like the forward–backward and Douglas–Rachford methods, as well as the more recent algorithm of Chambolle and Pock \cite{Chambolle.Pock2010}. 

If a minimizer to problem (\ref{eq:problem2}) exists, algorithm (\ref{eq:alg3}) converges for $0<\alpha<2/L$ and $0<\beta<1/\|A\|^2$ \cite{Loris.Verhoeven2011,Chen2013}. In the following section, we will prove convergence of algorithm (\ref{eq:alg4}) under the same conditions.

Algorithm~\ref{alg:npdaaverage} is very similar to the one used in \cite{Beck.Teboulle2009} for the special case of so-called Total Variation image denoising and deblurring problems. The main difference lies in the absence (in \cite{Beck.Teboulle2009}) of a feedback strategy for the inner loop: The authors of \cite{Beck.Teboulle2009} restart the inner iteration at $v_n^0=0$ (for all $n$) and observe that this, in combination with a fixed number of inner iterations, may lead to non-convergence of the outer loop. The main contribution of this paper therefore is the convergence resulting from the feedback strategy $v_{n}^0=v_{n-1}^{k_{\max}}$ (for all $n>0$).



\section{Convergence Results}\label{sec:results}

Three further lemma's are needed for proving convergence of algorithm (\ref{eq:alg4}).
\begin{lemma}\label{lemma:minimizer}
The minimizers $\hat u$ of problem (\ref{eq:problem2}) are characterized by the equations
\begin{equation}\label{eq:minimizers}
\left\{
\begin{array}{l}
\hat u=\prox_{\alpha h}(\hat u-\alpha \grad f(\hat u)- \alpha A^T\hat v)\\[3mm]
\hat v=\prox_{\beta\alpha^{-1} g^\ast}\left( \hat v+\beta\alpha^{-1} A\hat u\right)
\end{array}
\right.
\end{equation}
for any $\alpha,\beta>0$.
\end{lemma}
\proof The minimizers of (\ref{eq:problem2}) are characterized by the inclusion $0\in\partial (f+h+g\circ A)(\hat u)$, or
\begin{equation} \label{tmp55}
0=\grad f(\hat u)+\hat w+A^T \hat v\qquad\text{with}\quad \hat w\in\partial h(\hat u)\quad \text{and}\quad \hat v\in\partial g(A\hat u).
\end{equation}
The two latter inclusions can also be written as (see Lemma~\ref{theoremMoreau}):
\begin{displaymath}
\hat u=\prox_{\alpha h}(\hat u+\alpha \hat w)\quad\text{and}\quad \hat v=\prox_{\beta\alpha^{-1} g^\ast}(\hat v+\beta\alpha^{-1} A\hat u)
\end{displaymath}
where $\alpha,\beta>0$ are arbitrary parameter and $g^\ast$ is the Fenchel dual of $g$. One obtains equations (\ref{eq:minimizers}) by using the first equation of (\ref{tmp55}) to eliminate $\hat w$.\qed

\begin{lemma}\label{lemmagradfne} Let $f:\realr^d \rightarrow \realr$ be a convex function with Lipschitz continuous gradient (constant $L$). It follows that $L^{-1}\grad f$ is firmly non-expansive:
\begin{equation} \label{fne}
\|\grad f(u)-\grad f(v)\|_2^2\leq L \langle \grad f(u)-\grad f(v),u-v\rangle\qquad\forall u,v\in\realr^d.
\end{equation}
\end{lemma}
\proof See \cite[Part 2, Chapter X, Th. 4.2.2]{Hiriart-Urruty1993}.
\qed

\begin{lemma}\label{proxproperty1} Let $h:\realr^d\rightarrow\realr\cup\{+\infty\}$ be a convex, proper, lower semi-continuous function. The equality $x^+=\prox_h(x^-+\Delta)$ is equivalent to the inequality:
\begin{displaymath} 
\|x^+-x\|_2^2\leq\|x^--x\|_2^2-\|x^+-x^-\|_2^2+2\langle x^+-x,\Delta\rangle +2h(x)-2h(x^+)
\end{displaymath}
for all $x\in\realr^d$.
\end{lemma}
\proof  $x^+=\prox_h(x^-+\Delta)$
\begin{displaymath}
\begin{array}{l}
\displaystyle\Leftrightarrow\quad x^+=\arg\min_{x}\frac{1}{2}\|x-(x^-+\Delta)\|_2^2+h(x)\\[3mm]
\displaystyle\Leftrightarrow\quad 0\in x^+-x^--\Delta+\partial h(x^+)\\[3mm]
\displaystyle\Leftrightarrow\quad x^-+\Delta-x^+\in \partial h(x^+)\\[3mm]
\displaystyle\Leftrightarrow\quad h(x)\geq h(x^+)+\langle x^-+\Delta-x^+,x-x^+\rangle\qquad\forall x\\[3mm]
\displaystyle\Leftrightarrow\quad \|x^+-x\|_2^2\leq\|x^--x\|_2^2-\|x^+-x^-\|_2^2+2\langle x^+-x,\Delta\rangle +2h(x)-2h(x^+)
\end{array}
\end{displaymath}
\qed

We are now ready to state and prove the main theorems.
\begin{theorem}\label{theorem:convergence4}
Let $0<\alpha<2/L$, $0<\beta<1/\|A\|^2$ and $k_{\max}\in\naturaln_0$.
If the optimization problem (\ref{eq:problem2}) admits a solution, the nested primal-dual algorithm (\ref{eq:alg4}) will converge to one.
\end{theorem}
\proof Let $\hat u\in\arg\min_u f(u)+g(Au)+h(u)$, i.e. there exists $\hat v$ such that equations (\ref{eq:minimizers}) are satisfied.

We use Lemma~\ref{proxproperty1} on the definition of $u_{n}^{k_{\max}}$ in algorithm~(\ref{eq:alg4}):
\begin{equation}\label{eq:finalu1}
\begin{array}{lcl}
\|u_{n}^{k_{\max}}-\hat u\|_2^2&\leq&\|u_n-\hat u\|_2^2-\|u_{n}^{k_{\max}}-u_n\|_2^2+2\alpha h(\hat u)-2\alpha h(u_{n}^{k_{\max}})\\[3mm]
&&\qquad-2\alpha\langle u_{n}^{k_{\max}}-\hat u, \grad f(u_n)+ A^T v_{n}^{k_{\max}} \rangle,
\end{array}
\end{equation}
on the definition of $u_{n}^{k}$ in algorithm~(\ref{eq:alg4}):
\begin{equation}\label{eq:finalu2}
\begin{array}{lcl}
\|u_{n}^{k}-u_{n}^{k+1}\|_2^2&\leq&\|u_n-u_{n}^{k+1}\|_2^2-\|u_{n}^{k}-u_n\|_2^2+2\alpha h(u_{n}^{k+1})-2\alpha h(u_{n}^{k})\\[3mm]
&&\qquad-2\alpha\langle u_{n}^{k}-u_{n}^{k+1}, \grad f(u_n)+ A^T v_{n}^{k} \rangle
\end{array}
\end{equation}
for $k:0\ldots k_{\max}-1$, and on the first line of equations (\ref{eq:minimizers}):
\begin{equation}\label{eq:finalu3}
\begin{array}{lcl}
\|\hat u-u_{n}^{0}\|_2^2&\leq&\|\hat u-u_{n}^{0}\|_2^2 -\|\hat u-\hat u\|_2^2-2\alpha\langle \hat u-u_{n}^{0}, \grad f(\hat u)+ A^T \hat v \rangle\\[3mm]
&&\qquad+2\alpha h(u_{n}^{0})-2\alpha h(\hat u).
\end{array}
\end{equation}
Applying Lemma~\ref{proxproperty1} again to the definition of $u_n^k$ in algorithm (\ref{eq:alg4}) 
\begin{displaymath}
\begin{array}{lcl}
\|u_{n}^{k}-\hat u\|_2^2 &\leq&\|u_n-\hat u\|_2^2 -\|u_{n}^k-u_n\|_2^2-2\alpha\langle u_{n}^{k}-\hat u, \grad f(u_n)+ A^T v_{n}^k\rangle\\[3mm]
&&\qquad +2\alpha h(\hat u)-2\alpha h(u_{n}^{k})
\end{array}
\end{displaymath}
for $k:1\ldots k_{\max}-1$ and to the first line of equations (\ref{eq:minimizers})
\begin{displaymath}
\begin{array}{lcl}
\|\hat u-u_{n}^{k}\|_2^2 &\leq&\|\hat u-u_{n}^{k}\|_2^2 -\|\hat u-\hat u\|_2^2-2\alpha\langle \hat u-u_{n}^{k}, \grad f(\hat u) +A^T \hat v\rangle\\[3mm]
&&\qquad +2\alpha h(u_{n}^{k})-2\alpha h(\hat u)
\end{array}
\end{displaymath}
together  yields:
\begin{equation}\label{eq:finalu}
\begin{array}{lcl}
\|u_{n}^{k}-\hat u\|_2^2 &\leq&\|u_n-\hat u\|_2^2 -\|u_{n}^k-u_n\|_2^2-2\alpha \langle u_{n}^{k}-\hat u, \grad f(u_n)-\grad f(\hat u) \rangle\\[3mm]
&&\qquad-2\alpha \langle u_{n}^{k}-\hat u,A^T (v_{n}^k-\hat v)\rangle
\end{array}
\end{equation}
for $k:1\ldots k_{\max}-1$.

Finally, we apply Lemma~\ref{proxproperty1} to the definition of $v_{n}^{k+1}$ in algorithm (\ref{eq:alg4}):
\begin{displaymath}
\begin{array}{lcl}
\|v_{n}^{k+1}-\hat v\|_2^2 &\leq&\|v_{n}^k-\hat v\|_2^2 -\|v_{n}^{k+1}-v_{n}^k\|_2^2+2\beta \alpha^{-1}\langle v_{n}^{k+1}-\hat v,Au_{n}^k \rangle \\[3mm]
&&\qquad+2\beta \alpha^{-1}g^\ast(\hat v)-2\beta \alpha^{-1}g^\ast(v_{n}^{k+1})
\end{array}
\end{displaymath}
and to the second equation in system (\ref{eq:minimizers}):
\begin{displaymath}
\begin{array}{lcl}
\|\hat v-v_{n}^{k+1}\|_2^2 &\leq&\|\hat v-v_{n}^{k+1}\|_2^2 -\|\hat v-\hat v\|_2^2+2\beta \alpha^{-1}\langle \hat v-v_{n}^{k+1},A\hat u \rangle \\[3mm]
&&\qquad+2\beta \alpha^{-1}g^\ast(v_{n}^{k+1})-2\beta \alpha^{-1}g^\ast(\hat v)
\end{array}
\end{displaymath}
which together give:
\begin{equation}\label{eq:finalv}
\begin{array}{lcl}
\|v_{n}^{k+1}-\hat v\|_2^2 \leq\|v_{n}^k-\hat v\|_2^2 -\|v_{n}^{k+1}-v_{n}^k\|_2^2+2\beta \alpha^{-1}\langle v_{n}^{k+1}-\hat v,A(u_{n}^k-\hat u) \rangle
\end{array}
\end{equation}
for $k:0\ldots k_{\max}-1$.

By adding the inequalities (\ref{eq:finalu1}), (\ref{eq:finalu2}), (\ref{eq:finalu3}), (\ref{eq:finalu}) and (\ref{eq:finalv}), and after canceling some terms and rearranging the remaining inner products, one obtains:
\begin{equation}\label{temp81}
\begin{array}{lcl}
\sum_{k=1}^{k_{\max}} \beta \|u_{n}^{k}-\hat u\|_2^2 +\sum_{k=0}^{k_{\max}-1} \alpha^2\|v_{n}^{k+1}-\hat v\|_2^2\leq \sum_{k=0}^{k_{\max}-1}\beta\|u_n-\hat u\|_2^2 \\[3mm]
-\beta\|u_{n}^{k}-u_{n}^{k+1}\|_2^2
-\beta\|u_{n}^{k}-u_n\|_2^2
-2\alpha\beta\langle u_{n}^{k}-\hat u, \grad f(u_n)-\grad f(\hat u)\rangle\\[3mm]
+\alpha^2\|v_{n}^k-\hat v\|_2^2 -\alpha^2\|v_{n}^{k+1}-v_{n}^k\|_2^2+2\alpha\beta\langle u_{n}^k-u_{n}^{k+1},A^T(v_{n}^{k+1}-v_{n}^{k})\rangle
\end{array}
\end{equation}
(verification of the last term in inequality (\ref{temp81}) takes a few lines but is straightforward; the details are omitted).
On the second line one can use the following bound:
\begin{displaymath}
\begin{array}{lcl}
-\|u_{n}^{k}-u_n\|_2^2
-2\alpha\langle u_{n}^{k}-\hat u, \grad f(u_n)-\grad f(\hat u)\rangle\\[3mm]
=\alpha^2 \|\grad f(u_n)-\grad f(\hat u)\|_2^2
-2\alpha\langle u_n-\hat u,\grad f(u_n)-\grad f(\hat u)\rangle\\[3mm]
\qquad-\|u_{n}^{k}-u_n+\alpha (\grad f(u_n)-\grad f(\hat u))\|_2^2\\[3mm]
\leq \alpha(\alpha-2/L) \|\grad f(u_n)-\grad f(\hat u)\|_2^2
-\|u_{n}^{k}-u_n+\alpha (\grad f(u_n)-\grad f(\hat u))\|_2^2
\end{array}
\end{displaymath}
where we have used the fact that $L^{-1}\grad f$ is firmly non expansive (Lemma~\ref{lemmagradfne}), while the scalar products on the last line of inequality (\ref{temp81}) can be replaced by
\begin{equation}\label{eq:innerprodbound}
\begin{array}{lcl}
2\alpha\beta\langle u_{n}^k-u_{n}^{k+1}, A^T(v_{n}^{k+1}-v_{n}^{k})\rangle
&=&\beta\|u_{n}^k-u_{n}^{k+1}\|_2^2+\beta\alpha^2\|A^T(v_{n}^{k+1}-v_{n}^{k})\|_2^2\\[3mm]
&&-\beta\|u_{n}^k-u_{n}^{k+1}-\alpha A^T(v_{n}^{k+1}-v_{n}^{k})\|_2^2.
\end{array}
\end{equation}
Hence we can deduce that
\begin{equation}\label{tmp82}
\begin{array}{lcl}
\sum_{k=0}^{k_{\max}-1} \beta \|u_{n}^{k+1}-\hat u\|_2^2 +\alpha^2 \|v_{n}^{k+1}-\hat v\|_2^2\leq \sum_{k=0}^{k_{\max}-1}\beta\|u_n-\hat u\|_2^2+\alpha^2\|v_{n}^k-\hat v\|_2^2 \\[3mm]
+\beta\sum_{k=0}^{k_{\max}-1}\alpha(\alpha-2/L) \|\grad f(u_n)-\grad f(\hat u)\|_2^2
-\|u_{n}^{k}-u_n+\alpha (\grad f(u_n)-\grad f(\hat u))\|_2^2\\[3mm]
+\alpha^2\sum_{k=0}^{k_{\max}-1} -\|v_{n}^{k+1}-v_{n}^k\|_2^2+\beta\|A^T(v_{n}^{k+1}-v_{n}^{k})\|_2^2\\[3mm]
-\beta\sum_{k=0}^{k_{\max}-1}\|u_{n}^k-u_{n}^{k+1}
-\alpha A^T(v_{n}^{k+1}-v_{n}^{k})\|_2^2.
\end{array}
\end{equation}
Now we use the convexity of  $\|u_{n+1}-\hat u\|_2^2$ (as a function of $u_{n+1}$) and the last line of algorithm (\ref{eq:alg4}) to write:
\begin{displaymath}
\|u_{n+1}-\hat u\|_2^2\leq \gamma\sum_{k=0}^{k_{\max}-1} \|u_{n}^{k+1}-\hat u\|_2^2
\end{displaymath}
where $\gamma=1/k_{\max}$. Together with inequality (\ref{tmp82}), we finally find:
\begin{equation}\label{eq:tmp71}
\begin{array}{lcl}
\beta\|u_{n+1}-\hat u\|_2^2+\alpha^2\gamma \|v_{n+1}^0-\hat v\|_2^2 \leq\beta\|u_n-\hat u\|_2^2+\alpha^2\gamma \|v_{n}^{0}-\hat v\|_2^2\\[3mm] +\alpha \beta \gamma(\alpha-2/L)\|\grad f(u_n)-\grad f(\hat u)\|_2^2-\beta\gamma\sum_{k=0}^{k_{\max}-1}\|u_{n}^k-u_n+\alpha (\grad f(u_n)-\grad f(\hat u))\|_2^2\\[3mm]
-\alpha^2\gamma\sum_{k=0}^{k_{\max}-1} \|v_{n}^{k+1}-v_{n}^{k}\|_A^2
-\beta\gamma\sum_{k=0}^{k_{\max}-1}\|u_{n}^k-u_{n}^{k+1}-\alpha A^T(v_{n}^{k+1}-v_{n}^{k})\|_2^2
\end{array}
\end{equation}
where we have used the relation $v_{n+1}^0=v_{n}^{k_{\max}}$ and the norm $\|v\|_A^2=\|v\|_2^2-\beta\|A^Tv\|_2^2$ (it is a norm because  $0<\beta<1/\|A\|^2$).

Relation (\ref{eq:tmp71}) and assumption $0<\alpha<2/L$ implies that the sequence $(u_n,v_{n}^0)_{n\in\naturaln}$ is bounded. Hence a limit point exists: $(u_{n_j},v_{n_j}^0)\stackrel{j\to\infty}{\longrightarrow}(u^\dagger, v^\dagger)$. By summing inequalities (\ref{eq:tmp71}) from $n=0$ until $n=N$ one deduces also that:
\begin{displaymath}
\begin{array}{l}
\|\grad f(u_n)-\grad f(\hat u)\|_2^2 \stackrel{n\to\infty}{\longrightarrow}0,\\[3mm]
\|u_{n}^k-u_n+\alpha (\grad f(u_n)-\grad f(\hat u))\|_2^2
\stackrel{n\to\infty}{\longrightarrow}0\qquad k:0\ldots k_{\max}-1,\\[3mm]
\|v_{n}^{k+1}-v_{n}^{k}\|_A^2\stackrel{n\to\infty}{\longrightarrow}0\qquad k:0\ldots k_{\max}-1,\\[3mm]
\|u_{n}^k-u_{n}^{k+1}-\alpha A^T(v_{n}^{k+1}-v_{n}^{k})\|_2^2\stackrel{n\to\infty}{\longrightarrow}0\qquad k:0\ldots k_{\max}-1.
\end{array}
\end{displaymath}
This in turn implies that 
\begin{displaymath}
u_{n_j}^{k+1}\stackrel{j\to\infty}{\longrightarrow}u^\dagger, \qquad
v_{n_j}^{k+1}\stackrel{j\to\infty}{\longrightarrow}v^\dagger,\qquad k:0\ldots k_{\max}-1
\end{displaymath}
also. It follows from the continuity of the operations in the right hand sides of algorithm (\ref{eq:alg4}) that $(u^\dagger, v^\dagger)$ satisfies the equations (\ref{eq:minimizers}), which characterize the minimizers of problem (\ref{eq:problem2}). One can then replace $(\hat u,\hat v)$ by $(u^\dagger, v^\dagger)$ in inequality (\ref{eq:tmp71}) to obtain
\begin{displaymath}
\beta\|u_{n+1}-u^\dagger\|_2^2+\alpha^2\gamma \|v_{n+1}^0-v^\dagger\|_2^2 \leq\beta\|u_n-u^\dagger\|_2^2+\alpha^2\gamma \|v_{n}^{0}-v^\dagger\|_2^2
\end{displaymath}
This then implies the convergence of the whole sequence $(u_n,v_n^0)$ to $(u^\dagger, v^\dagger)$.\qed

As usual, one expects that better convergence results can be obtained when strong convexity of the objective function is assumed. 
In \cite[Example 27.12]{Bauschke2011} the linear convergence rate (to the unique minimizer $\hat u$ of problem (\ref{eq:problem1})) of the proximal-gradient algorithm (\ref{eq:proxgrad}) is proven, when  $f$ is strongly convex (parameter $\mu$), $\grad f$ is Lipschitz continuous (parameter $L$) and $0<\alpha<2/L$:
\begin{displaymath}
\|u_{n+1}-\hat u\|_2^2\leq(1+\mu \alpha(\alpha L-2)) \|u_n-\hat u\|_2^2,
\end{displaymath}
where $0\leq 1+\mu \alpha(\alpha L-2)<1$.
In \cite{Chaux.Pesquet.ea2009} linear convergence rate of the proximal gradient algorithm (\ref{eq:proxgrad}) is shown for strongly convex $h$ (instead of $f$). The following theorem thus complements the results of \cite{Loris.Verhoeven2011,Chen2013,Chen2016,Chaux.Pesquet.ea2009}.

\begin{theorem}\label{theorem:linearconvergence2} 
Let $0<\alpha<2/L$, $0<\beta<1/\|A\|^2$ and $k_{\max}\in\naturaln_0$.
In addition, we assume that $f$ is strongly convex (parameter $\mu$), that $h=0$ and that $A^T $ is coercive (parameter $\sigma>0$): $\|A^Tv\|_2\geq \sigma\|v\|_2$ for all $v \in\realr^{d'}$.

Then the primal-dual algorithm (\ref{eq:alg4}) converges to the minimizer $\hat u$ of problem (\ref{eq:problem2}) at a linear linear rate:
\begin{displaymath}
\beta\|u_{n+1}-\hat u\|_2^2+\alpha^2\gamma \|v_{n+1}^0-\hat v\|_2^2 \leq\epsilon \left(\beta\|u_n-\hat u\|_2^2+\alpha^2\gamma \|v_{n}^{0}-\hat v\|_2^2\right)
\end{displaymath}
for some $0\leq \epsilon<1$ and $\gamma=1/k_{\max}$. Here, $\hat v$ is the dual variable of equations (\ref{eq:minimizers}).
\end{theorem}
\proof As $f$ is strongly convex, problem (\ref{eq:problem2}) is guaranteed to have a (unique) solution $\hat u$. Hence, there exists $\hat v$ such that equations (\ref{eq:minimizers}) are satisfied.

We start from inequality (\ref{temp81}) derived in the proof of Theorem~\ref{theorem:convergence4}.
We use the following bound:
\begin{displaymath}
\begin{array}{l}
\|u_n-\hat u\|_2^2-\|u_{n}^{k}-u_n\|_2^2+2\alpha\langle u_{n}^{k}-\hat u, \grad f(\hat u)-\grad f(u_n)\rangle\\[3mm]
 =  \|u_n-\hat u\|_2^2+\alpha^2 \|\grad f(u_n)-\grad f(\hat u)\|_2^2-2\alpha \langle u_n-\hat u, \grad f(u_n)-\grad f(\hat u)\rangle\\[3mm]
\qquad\qquad\qquad\qquad -\|u_{n}^{k}-u_n+\alpha(\grad f(u_n)-\grad f(\hat u))\|_2^2\\[3mm]
\leq \|u_n-\hat u\|_2^2+(\alpha^2 L-2\alpha) \langle u_n-\hat u, \grad f(u_n)-\grad f(\hat u)\rangle\\[3mm]
\qquad\qquad\qquad\qquad-\|u_{n}^{k}-u_n+\alpha(\grad f(u_n)-\grad f(\hat u))\|_2^2\\[3mm]
\leq (1+\mu \alpha(\alpha L-2)) \|u_n-\hat u\|_2^2-\|u_{n}^{k}-u_n+\alpha(\grad f(u_n)-\grad f(\hat u))\|_2^2
\end{array}
\end{displaymath}
where we have used the fact that $L^{-1}\grad f$ is firmly non expansive (lemma~\ref{lemmagradfne}), 
the relation
\begin{displaymath}
\langle \grad f(u_n)-\grad f(\hat u),u_n-\hat u \rangle\geq \mu \|u_n-\hat u\|_2^2,
\end{displaymath}
which is a consequence of the strong convexity of $f$:
\begin{displaymath}
\begin{array}{lcl}
f(u_n)&\geq& f(\hat u)+\langle \grad f(\hat u),u_n-\hat u \rangle+\mu \|u_n-\hat u\|_2^2 /2\\[3mm]
f(\hat u)&\geq& f(u_n)+\langle \grad f(u_n),\hat u-u_n \rangle+\mu \|u_n-\hat u\|_2^2 /2,
\end{array}
\end{displaymath} 
and the assumption that $\alpha L-2<0$. The second term on the right hand side can be written as:
\begin{displaymath}
\|u_{n}^{k}-u_n+\alpha(\grad f(u_n)-\grad f(\hat u))\|_2^2=\alpha^2 \|A^T(v_{n}^k-\hat v)\|_2^2
\end{displaymath}
on account of the assumption that $h=0$ (in this case $\prox_{\alpha h}$ is the identity), the definition of $u_n^k$  in (\ref{eq:alg4}) and the first line in equations (\ref{eq:minimizers}).

We thus find from inequality (\ref{temp81}):
\begin{displaymath}
\begin{array}{l}
 \sum_{k=0}^{k_{\max}-1} \beta\|u_{n}^{k+1}-\hat u\|_2^2 + \alpha^2\|v_{n}^{k+1}-\hat v\|_2^2\\[3mm]
 \leq \sum_{k=0}^{k_{\max}-1}\beta(1+\mu \alpha(\alpha L-2))\|u_n-\hat u\|_2^2
-\alpha^2\beta \|A^T(v_{n}^k-\hat v)\|_2^2-\beta\|u_{n}^{k}-u_{n}^{k+1}\|_2^2\\[3mm]
+\alpha^2\|v_{n}^k-\hat v\|_2^2 -\alpha^2\|v_{n}^{k+1}-v_{n}^k\|_2^2
+2\alpha\beta\langle u_{n}^k-u_{n}^{k+1},A^T(  v_{n}^{k+1}-v_{n}^{k})\rangle.
\end{array}
\end{displaymath}
The inner products on the last line can bounded using the first two terms in the right hand side of expression (\ref{eq:innerprodbound}), such that one finds:
\begin{displaymath}
\begin{array}{l}
 \sum_{k=0}^{k_{\max}-1} \beta\|u_{n}^{k+1}-\hat u\|_2^2 + \alpha^2\|v_{n}^{k+1}-\hat v\|_2^2\\[3mm]
\leq \sum_{k=0}^{k_{\max}-1}\beta(1+\mu \alpha(\alpha L-2))\|u_n-\hat u\|_2^2
+\alpha^2\|v_{n}^k-\hat v\|_A^2 -\alpha^2\|v_{n}^{k+1}-v_{n}^k\|_A^2
\end{array}
\end{displaymath}
using the norm $\|v\|_A^2=\|v\|_2^2-\beta\|A^Tv\|_2^2$. The last term on the right hand side is dropped and the definition $u_{n+1}=\gamma\sum_{k=0}^{k_{\max}-1} u_{n}^{k+1}$ (with $\gamma=k_{\max}^{-1}$) and the convexity of $\|\cdot-\hat u\|_2^2$ then imply that:
\begin{displaymath}
\begin{array}{l}
  \beta\|u_{n+1}-\hat u\|_2^2 + \sum_{k=0}^{k_{\max}-1}\alpha^2\gamma\|v_{n}^{k+1}-\hat v\|_2^2\leq \\[3mm]
  \qquad\qquad\qquad \beta(1+\mu \alpha(\alpha L-2))\|u_n-\hat u\|_2^2   +\sum_{k=0}^{k_{\max}-1}\alpha^2\gamma\|v_{n}^k-\hat v\|_A^2.
\end{array}
\end{displaymath}
As $A^T$ is coercive, one has that $\|v_{n}^k-\hat v\|_A^2\leq (1-\beta \sigma^2)\|v_{n}^k-\hat v\|_2^2$ such that:
\begin{displaymath}
\begin{array}{l}
\beta\|u_{n+1}-\hat u\|_2^2+\alpha^2\gamma \|v_{n+1}^0-\hat v\|_2^2 \leq\\[3mm]
\qquad\qquad\qquad \beta(1+\mu \alpha(\alpha L-2))\|u_n-\hat u\|_2^2+\alpha^2\gamma(1-\beta \sigma^2) \|v_{n}^{0}-\hat v\|_2^2

\end{array}
\end{displaymath}
as $v_{n+1}^0=v_n^{k_{\max}}$.
By setting $\epsilon=\max((1+\mu \alpha(\alpha L-2)), 1-\beta \sigma^2)$ one finds the announced inequality.

In order to show that $0\leq \epsilon <1$, one proceeds as follows. 
On the one hand, $\alpha L-2<0$ implies $1+\mu \alpha(\alpha L-2)<1$, while $1+\mu \alpha(\alpha L-2)$ reaches a minimum for $\alpha=1/L$. This minimum is $1-\mu/L\geq 0$ as $\mu\leq L$ for strongly convex functions $f$ (with parameter $\mu$) with Lipschitz continuous gradient (parameter $L$). One sees that $0<1-\beta \sigma^2<1$ on account of $0<\sigma^{2}\leq\|A\|^{2}<\beta^{-1}$.
\qed

The proof of Theorem~\ref{theorem:linearconvergence2} unfortunately requires the assumption that $h=0$. When $h=0$ the dual problem (\ref{eq:dualproblem}) reduces to a quadratic plus proximable term:
\begin{displaymath}
\min_v \frac{1}{2}\|a/\alpha-A^T v\|_2^2 +\frac{1}{\alpha}g^\ast (v).
\end{displaymath}
When $A^T$ is coercive, the first term is strongly convex. For general $h$ we conjecture that a linear convergence rate still holds for algorithm (\ref{eq:alg4}) when one assumes, in addition to the strong convexity of $f$, that the function $\psi$ appearing in the dual problem (\ref{eq:dualproblem}) is strongly convex. In \cite{Chen2016} a linear convergence rate is shown for algorithm (\ref{eq:alg4}) with $k_{\max}=1$ and assuming that $g^\ast$ is strongly convex (in addition to some other assumptions on $f$ and $A$).

\section{Conclusions}

A generalization of the proximal gradient algorithm (\ref{eq:proxgrad}) consisting of nested primal and dual iterations was discussed and convergence was shown. The iterative algorithm requires access to a gradient and two proximal operators, but not to the inverse of the linear operator appearing in problem (\ref{eq:problem2}). Similar problems and related algorithms are also discussed in \cite{Zhang.Burger.ea2011,Combettes2014}. Under some additional conditions (related to strong convexity of the cost function) a linear convergence rate was shown.

Nested iterative algorithms are abundant in numerical and applied mathematics. The proposed algorithm is very similar to the one in \cite{Beck.Teboulle2009}. The main novelty lies in the rigorous discussion of the inner loop starting and stopping criterion. One often encounters inner loop stopping criteria of the form $\|v_n^k-v_n^{k+1}\|_2<\epsilon$, which may give satisfactory numerical results, but may not guarantee convergence of the outer loop. In addition, such a condition may just indicate slow convergence of the inner loop. Additionally, the inner loop starting point is often neglected in theoretical descriptions. In practice (i.e. in the implementation code instead of in papers), a feedback/warm start mechanism  of type $v_{n}^0=v_{n-1}^{k_{\max}}$ is sometimes added to ``speed up" convergence of the inner loop. Here we have shown that such a small change can already be sufficient to guarantee convergence. Such a discussion has not been given before, and it is the main contribution of this paper: In the context of nested iterative algorithms, inner loop ``starting rules" can be just as effective as ``stopping rules'' for guaranteeing convergence.

No numerical experiments are presented. In fact, the proposed Algorithm~\ref{alg:npdaaverage} cannot be expected to be  state-of-the-art by itself (lack of variable step length or line-search strategies \cite{Chouzenoux-etal-2014,Bonettini2009}). The point here is just to prove that the described mechanism is sufficient for convergence. More sophisticated algorithms, e.g. incorporating line-search rules to speed-up convergence, exist. In those cases too, one could investigate the role of the warm start mechanism on the convergence of nested iterations.

Another possible extension concerns the convergence of nested \emph{accelerated} primal-dual algorithms of Nesterov type. The convergence of the iterates of accelerated projected-gradient \cite{Nesterov1983a} and proximal-gradient algorithms \cite{Beck.Teboulle2008} was shown in \cite{Chambolle2015}. Nested primal-dual versions were proposed in \cite{Beck.Teboulle2009}, again without feedback. The proof of convergence of algorithms of that type is still an open problem. 

Another generalization concerns the use of variable stepsizes ($\alpha_n$ instead of $\alpha$ and $\beta_n$ instead of $\beta$) in algorithm (\ref{eq:alg4}). Finally, the condition $0<\beta<1/\|A\|^2$ seems to be too restrictive in view of the step size condition in Lemma~\ref{lemma:dual} ($0<\beta<2/\|A\|^2$). A variation of algorithm (\ref{eq:alg4}) with a different feedback strategy will be described in \cite{Chen2018}.

\section{Acknowledgements}
JC is sponsored by the China Scholarship Council. IL is a Research Associate of the Fonds de la Recherche Scientifique - FNRS and is also supported by a ULB ARC grant.

\bibliography{biblio}{}%

\begin{thebibliography}{10}

\bibitem{Goldstein1964}
A.~A. Goldstein.
\newblock Convex programming in {Hilbert} space.
\newblock {\em Bulletin of the American Mathematical Society}, 70:709--710,
  1964.

\bibitem{ComWa2005}
Patrick~L. Combettes and Valerie~R. Wajs.
\newblock Signal recovery by proximal forward-backward splitting.
\newblock {\em Multiscale Model. Simul.}, 4(4):1168--1200, January 2005.

\bibitem{Moreau1965}
J.~J. Moreau.
\newblock Proximité et dualité dans un espace hilbertien.
\newblock {\em Bull. Soc. Math. France}, 93:273--299, 1965.

\bibitem{Combettes2014a}
Patrick~L. Combettes and Băng~C. Vũ.
\newblock Variable metric forward–backward splitting with applications to
  monotone inclusions in duality.
\newblock {\em Optimization}, 63(9):1289--1318, 2014.

\bibitem{Chouzenoux-etal-2014}
E.~Chouzenoux, J.-C. Pesquet, and A.~Repetti.
\newblock Variable metric forward-backward algorithm for minimizing the sum of
  a differentiable function and a convex function.
\newblock {\em J. Optim. Theory Appl.}, 162(1):107--132, July 2014.

\bibitem{Condat2013}
L.~Condat.
\newblock A primal-dual splitting method for convex optimization involving
  {L}ipschitzian, proximable and linear composite terms.
\newblock {\em J. Optimization Theory and Applications}, 158(2):460--479, 2013.

\bibitem{Combettes2014b}
P.~L. Combettes, L.~Condat, J.-C. Pesquet, and B.~C. Vu.
\newblock A forward-backward view of some primal-dual optimization methods in
  image recovery.
\newblock In {\em 2014 IEEE International Conference on Image Processing
  (ICIP)}, pages 4141--4145, 2014.

\bibitem{Zhang.Burger.ea2011}
Xiaoqun Zhang, Martin Burger, and Stanley Osher.
\newblock A unified primal-dual algorithm framework based on bregman iteration.
\newblock {\em J Sci Comput}, 46:20--46, 2011.

\bibitem{Combettes2012}
Patrick~L. Combettes and Jean-Christophe Pesquet.
\newblock Primal-dual splitting algorithm for solving inclusions with mixtures
  of composite, lipschitzian, and parallel-sum type monotone operators.
\newblock {\em Set-Valued and Variational Analysis}, 20(2):307--330, 2012.

\bibitem{Chen2016}
Peijun Chen, Jianguo Huang, and Xiaoqun Zhang.
\newblock A primal-dual fixed point algorithm for minimization of the sum of
  three convex separable functions.
\newblock {\em Fixed Point Theory and Applications}, 2016(1):54, 2016.

\bibitem{Chambolle2004}
Antonin Chambolle.
\newblock An algorithm for total variation minimization and applications.
\newblock {\em Journal of Mathematical Imaging and Vision}, 20:89--97, 2004.

\bibitem{Chambolle2005}
A.~Chambolle.
\newblock Total variation minimization and a class of binary mrf models.
\newblock In {\em Energy Minimization Methods in Computer Vision and Pattern
  Recognition}, volume 3757 of {\em Lecture Notes in Computer Science}, pages
  136--152, 2005.

\bibitem{Aujol2009}
Jean-François Aujol.
\newblock Some first-order algorithms for total variation based image
  restoration.
\newblock {\em J Math Imaging Vis}, 34:307--327, 2009.

\bibitem{Beck.Teboulle2009}
A.~Beck and M.~Teboulle.
\newblock Fast gradient-based algorithms for constrained total variation image
  denoising and deblurring problems.
\newblock {\em Image Processing, IEEE Transactions on}, 18(11):2419 --2434,
  nov. 2009.

\bibitem{Bonettini2015}
S.~Bonettini, I.~Loris, F.~Porta, and M.~Prato.
\newblock Variable metric inexact line-search based methods for nonsmooth
  optimization.
\newblock {\em Siam Journal on Optimization}, 26(2):891--921, 2016.

\bibitem{Bonettini2016}
S.~Bonettini, I.~Loris, F.~Porta, M.~Prato, and S.~Rebegoldi.
\newblock On the convergence of a linesearch based proximal-gradient method for
  nonconvex optimization.
\newblock {\em Inverse Problems}, 33(5):055005, 2017.

\bibitem{Salzo2012}
Saverio Salzo and Silvia Villa.
\newblock Inexact and accelerated proximal point algorithms.
\newblock {\em Journal of Convex Analysis}, 19(4):1167--1192, 2012.

\bibitem{Schmidt2011}
Mark Schmidt, Nicolas~Le Roux, and Francis Bach.
\newblock Convergence rates of inexact proximal-gradient methods for convex
  optimization.
\newblock In {\em Proceedings of the 24th International Conference on Neural
  Information Processing Systems}, NIPS'11, pages 1458--1466, USA, 2011. Curran
  Associates Inc.

\bibitem{Hu2016}
Yue Hu, Eric~C. Chi, and Genevera~I. Allen.
\newblock {\em Splitting Methods in Communication, Imaging, Science, and
  Engineering}, chapter ADMM Algorithmic Regularization Paths for Sparse
  Statistical Machine Learning, pages 433--460.
\newblock Springer, 2016.

\bibitem{Rose2015}
S.~Rose, MS~Andersen, EY~Sidky, and X~Pan.
\newblock Noise properties of {CT} images reconstructed by use of constrained
  total-variation, data-discrepancy minimization.
\newblock {\em Medical Physics}, 42(5):2690--2698, 2015.

\bibitem{Rose2016}
Sean Rose, Martin~S. Andersen, Emil~Y. Sidky, and Xiaochuan Pan.
\newblock Technical note: Proximal ordered subsets algorithms for {TV}
  constrained optimization in {CT} image reconstruction.
\newblock Technical report, The University of Chicago, 2016.
\newblock arXiv:1603.08889v1.

\bibitem{Loris.Verhoeven2011}
Ignace Loris and Caroline Verhoeven.
\newblock On a generalization of the iterative soft-thresholding algorithm for
  the case of non-separable penalty.
\newblock {\em Inverse Problems}, 27(12):125007, 2011.

\bibitem{Chen2013}
Peijun Chen, Jianguo Huang, and Xiaoqun Zhang.
\newblock A primal-dual fixed point algorithm for convex separable minimization
  with applications to image restoration.
\newblock {\em Inverse Problems}, 29(2):025011--, 2013.

\bibitem{Bauschke2011}
Heinz~H. Bauschke and Patrick~L. Combettes.
\newblock {\em Convex Analysis and Monotone Operator Theory in Hilbert Spaces}.
\newblock CMS book in mathematics. Springer, 2011.

\bibitem{Combettes2014}
P.~L. Combettes, L.~Condat, J.-C. Pesquet, and B.~C. Vu.
\newblock A forward-backward view of some primal-dual optimization methods in
  image recovery.
\newblock In {\em 2014 IEEE International Conference on Image Processing
  (ICIP)}, pages 4141--4145, 2014.

\bibitem{Condat2014}
L.~Condat.
\newblock A generic proximal algorithm for convex optimization – application
  to total variation minimization.
\newblock {\em IEEE Signal Proc. Letters}, 21(8):1054--1057, 2014.

\bibitem{Chambolle.Pock2010}
Antonin Chambolle and Thomas Pock.
\newblock A first-order primal-dual algorithm for convex problems with
  applications to imaging.
\newblock {\em J Math Imaging Vis}, 40:120--145, 2011.

\bibitem{Hiriart-Urruty1993}
J.~B. Hiriart-Urruty and C.~Lemarechal.
\newblock {\em Convex analysis and minimization algorithms}.
\newblock Springer, 1993.

\bibitem{Chaux.Pesquet.ea2009}
Caroline Chaux, Jean-Christophe Pesquet, and Nelly Pustelnik.
\newblock Nested iterative algorithms for convex constrained image recovery
  problems.
\newblock {\em SIAM J. Imaging Sci.}, 2(2):730--762, January 2009.

\bibitem{Bonettini2009}
S.~Bonettini, R.~Zanella, and L.~Zanni.
\newblock A scaled gradient projection method for constrained image deblurring.
\newblock {\em Inverse Problems}, 25(1):015002, 2009.

\bibitem{Nesterov1983a}
Yu~E. Nesterov.
\newblock A method for solving a convex programming problem with convergence
  rate {$\mathcal{O}(1/k^2)$}.
\newblock {\em Soviet Math. Dokl.}, 27:372--376, 1983.

\bibitem{Beck.Teboulle2008}
Amir Beck and Marc Teboulle.
\newblock A fast iterative shrinkage-threshold algorithm for linear inverse
  problems.
\newblock {\em SIAM Journal on Imaging Sciences}, 2:183--202, 2009.

\bibitem{Chambolle2015}
A.~Chambolle and Ch. Dossal.
\newblock On the convergence of the iterates of the ``fast iterative
  shrinkage/thresholding algorithm".
\newblock {\em Journal of Optimization Theory and Applications},
  166(3):968--982, 2015.

\bibitem{Chen2018}
Jixin Chen.
\newblock {\em Domain decomposition methods and convex optimization with
  applications to inverse problems}.
\newblock PhD thesis, East China Normal University and Université libre de
  Bruxelles, 2018.
\newblock In preparation.

\end{thebibliography}
\bibliographystyle{unsrt}

\end{document}